\documentclass{article}

\usepackage{arxiv}

\usepackage{hyperref} 
\usepackage{url} 
\usepackage{booktabs} 
\usepackage{amsfonts} 
\usepackage{nicefrac} 
\usepackage{microtype} 
\usepackage{lipsum}

\usepackage{graphicx}
\usepackage{amsmath}

\usepackage{xcolor}

 \usepackage[T2A]{fontenc}
\usepackage[english]{babel}

\title{Fast sixth-order algorithm based on the generalized Cayley transform for the Zakharov-Shabat system in optical applications}

\author{
 Sergey Medvedev$^{1,2,*}$, Igor Chekhovskoy$^{2}$, Irina Vaseva$^{1,2}$, Mikhail Fedoruk$^{2,1}$\\
$^{1}$ Federal Research Center for Information and Computational Technologies,\\ Novosibirsk
630090, Russia,\\
$^{2}$ Novosibirsk State University, Novosibirsk 630090, Russia,\\
* Corresponding author: medvedev@ict.nsc.ru
}

\begin{document}
\maketitle

\begin{abstract}
Based on the generalized Cayley transform, a family of conservative one-step schemes of the sixth order of accuracy for the Zakharov-Shabat system is constructed. The exponential integrator is a special case. Schemes based on rational approximation allow the use of fast algorithms to solve the initial problem for a large number of values of the spectral parameter.
\end{abstract}


\section{Introduction}
In quantum mechanics and optics, problems arise that are described by linear ordinary differential equations with variable coefficients. The solutions of such equations can be written analytically only in exceptional cases. Therefore, to solve even linear ordinary differential equations, one has to use numerical methods, which are offered in a large number and variety (see, for example, fundamental books~\cite{hairer1993n, hairer1996solving, hairer2006geometric} and bibliography there).

Our work is also devoted to the construction of numerical methods for solving systems of linear ordinary differential equations with variable coefficients. The systems under consideration have features that require the construction of specialized numerical methods with the ability to perform fast massive computations. The main system and the problem for which we are developing our methods is to solve the direct spectral problem for the Zakharov-Shabat system (ZS).

The great interest in solving the direct spectral problem for the ZS system is based on the fact that it is the first step of the inverse scattering problem method for solving the nonlinear Schr\"odinger equation (NLSE) and its integrable modifications~\cite{hasegawa1995solitons}. The idea to use solitons for data transmission in optical fiber lines arose for the first time in~\cite{hasegawa1973transmission}. After this work, the nonlinear Schr\"odinger equation and its modifications were extremely intensively studied in relation to fiber telecommunication systems \cite{hasegawa2003optical,mollenauer2006solitons,turitsyn2003physics,yushko2014timing}.
Later, the idea was put forward to use multisoliton pulses in fiber-optic data transmission lines, when information is modulated and restored in the so-called nonlinear Fourier domain \cite{hasegawa1993eigenvalue,yousefi2014information,hari2016multieigenvalue}. Despite the fact that NLSE is an integrable system \cite{shabat1972exact}, its numerous studies were carried out by numerical methods. A classic overview of numerical methods for NLSE is given in \cite{taha1984analytical}.
The next step for the study of the direct spectral problem was made in the papers \cite{boffetta1992computation, burtsev1998numerical}, which were devoted to the numerical determination of scattering data for the ZS system.

It should be emphasized that the ZS system appears in other optical applications \cite{lamb1980elements}. In particular, the problem of scattering by Bragg gratings, which serve as the basis for optical filters in high-speed fiber data transmission lines, is reduced to the ZS system \cite{kashyap1999fibre,podivilov2006exactly}. The Schr\"odinger equation for two-level quantum systems with a time-dependent Hamiltonian takes the form of the ZS system \cite{akulin2005coherent,carmel2000geometrical}.
Moreover, the first example corresponds to the normal dispersion in the NLSE, and the second example corresponds to the anomalous dispersion. 
 In both cases, the ZS system has the skew-gradient form and preserves the quadratic integral \cite{medvedev2019exponential,medvedev2020exponential}. This integral is positively defined for the second case. 
In addition, the nonlinear Fourier transform is used to analyze coherent structures in dissipative systems and laser radiation~\cite{chekhovskoy2019nonlinear,sugavanam2019analysis}.

At present, the main goal of numerical methods for the direct spectral problem is the construction of fast methods, which, apparently, were first proposed in \cite{wahls2013introducing,wahls2015fast}. The current state of numerical methods and the prospects for the application of fast nonlinear Fourier transform for data transmission are given in \cite{turitsyn2017nonlinear,vasylchenkova2019direct,span2019time}.
The basic idea of \cite{wahls2013introducing, wahls2015fast} and subsequent works in this direction is to reduce the transition matrix or the matrix of the fundamental solution to a polynomial in the spectral parameter with matrix coefficients. Then the calculation of the spectral data for the continuous spectrum is reduced to the calculation of this polynomial for a large number of points of the continuous spectrum, and the calculation of the discrete spectrum is reduced to the calculation of the roots of the resulting polynomial. The advantage of this approach is the ability to use fast algorithms to compute polynomials \cite{pan2012structured,mcnamee2013numerical}.

For effective and practical application of the nonlinear Fourier transform for data transmission in fiber lines, high-precision and fast methods for solving the direct spectral problem for the ZS system are required. At the moment, the authors are aware of several fourth-order schemes \cite{chimmalgi2019fast,vaibhav2019efficient,medvedev2019exponential,medvedev2020conservative}, 
which allow the use of fast algorithms for multipoint computation of polynomials and finding the roots. Moreover, the scheme used in \cite{chimmalgi2019fast} is applied on an irregular grid, in \cite{vaibhav2019efficient} the Runge-Kutta method is used, which also requires the calculation of values within each elementary cell of the grid. Therefore, the above schemes require either interpolation within a unit cell or computation on a grid with a large step size. At the same time, a feature of the problem being solved is that the ZS system is specified in a tabular form on an equidistant grid. Numerical experiments in \cite{chimmalgi2019fast} showed unsatisfactory results when interpolating over several neighboring points. And only global interpolation gives satisfactory results. Our task is to construct schemes on an equidistant grid without using any interpolation. There are also schemes of the sixth \cite{blanes2006fourth,blanes2017high} and the eighth order \cite{alvermann2011high} for solving the general non-autonomous system. However, sixth-order schemes violate the unitarity of the transition matrix and, as a result, do not preserve the quadratic integral. Eighth-order schemes are constructed for expansion in Legendre polynomials and, therefore, require knowledge of the expansion coefficients. In \cite{chimmalgi2019fast}, Richardson interpolation was used for the 4th order scheme from \cite{blanes2017high}, which made it possible to build a fast 6th order algorithm. 

In our work, we consider the linear system
\begin{equation}\label{Psit}
\frac{d\Psi}{dt}=Q(t)\Psi,
\end{equation}
where the matrix $Q(t)$ depends on time $t$. Such systems arise in many optical and physical applications, as described above. The main attention will be paid to the direct spectral problem of the ZS system.

Our idea for constructing one-step methods 
\begin{equation}
\Psi_{n+1}=T\Psi_{n},    
\end{equation}
where $\Psi_n=\Psi(n\tau)$,
consists of finding first the expansion of the transition matrix $T$ into a Maclaurin series in terms of a small parameter, which is taken as the grid step size $\tau$, with the required accuracy, and then consistently replacing the derivatives with different analogs \cite{medvedev2019exponential,medvedev2020exponential}. Using this expansion, it is possible to determine what order of the derivatives of the matrix $Q(t)$ 
is necessary to construct a difference scheme of a given order of accuracy. By adding higher-order terms, one can represent the transition matrix in different forms. In particular, in exponential form, in the form of a product of exponentials, or a more general exponential expansion. Further, these functions can be approximated through rational functions. This work is devoted to the construction of a rational approximation of the transition matrix. For the exact conservation of quadratic integrals by schemes, it is proposed to use an approximation in the form of a generalized Cayley transform.

Let's formulate the main features of the direct spectral problem of the ZS system:

1. The matrix $Q(t)$ of the system (\ref{Psit}) is given on a uniform grid with a step size $\tau$, so the problem arises of constructing difference schemes that use only the values of $Q(t)$ at the grid nodes. If the values of $Q(t)$ can be calculated at any point $t$, then it is reasonable to calculate them within the integration step. This is exactly what is done when using Runge-Kutta schemes. If exponential integrators based on the Magnus expansion are used, then to approximate the integrals in this expansion, quadrature formulas are used at the optimal nodes within the grid cell.

2. The matrix $Q(t)$ is a polynomial in the complex parameter $\zeta$ and it is required to integrate the equation for a large number of values of $\zeta$. Therefore, it becomes necessary to represent the product of transition matrices as a polynomial with matrix coefficients and use fast algorithms to calculate them for a large number of values. This dictates the choice of a special kind of transition matrices.

3. For real spectral parameters $\zeta$ the ZS system conserves the quadratic integral, therefore the scheme must satisfy this property as well. Since the integration takes place over a large area, this conservation must be accurate.
This class of equations includes the Schr\"odinger equation with a time-dependent Hamiltonian.

4. The solutions of the system for the spectral parameters $\zeta$ lying in the upper complex half-plane have exponentially increasing and decaying solutions, therefore the scheme must be A-stable and, according to the condition of the second Dahlquist barrier, this condition is satisfied only by one-step explicit schemes.

5. System ZS and AKNS have dimension two. This allows constructing schemes containing exponents and other functions from matrices without significantly increasing the computational cost.

The article is organized as follows. Section~\ref{gen} contains general schematics for the system~(\ref{Psit}). In section~\ref{ZS}, the general theory is applied to the ZS system. Numerical experiments for the ZS system are given in section~\ref{num}.

\section{General theory of schemes}\label{gen}

Exponential difference schemes are based on the Magnus expansion \cite{iserles2000lie,blanes2009magnus}. The Magnus expansion contains integration over a time interval \cite{magnus1954exponential}; therefore, all schemes are based on the approximation of multiple integrals using cubature formulas on a set of nodes within an elementary subinterval. However, another option is also possible, in which the integrand is replaced by an expansion in a Taylor series, and the subsequent integration of this expansion is performed. If the derivatives from the Taylor series are approximated by difference analogs with suitable accuracy, then we obtain a difference scheme. For difference analogs, one can use the values of the system matrix $Q(t)$ only for $t$, in which the matrix $Q$ is given. This allows one to explicitly exclude interpolation within an elementary subinterval. This approach was used in \cite{medvedev2019exponential,medvedev2020exponential} and will be applied in this work.

The Magnus expansion transforms into an exact exponential solution for a system with a constant matrix. However, calculating the matrix exponential requires significant computational resources for high-dimensional matrices \cite{moler2003nineteen,najfeld1995derivatives}. Therefore, the idea arose to use rational approximations for systems with constant coefficients \cite{blue1970rational}. Among rational approximations, the diagonal Pad\'e approximation stands out \cite{calahan1967numerical}. Recently, the idea of using diagonal Pad\'e approximation was applied to the linear Schr\"odinger equation with the time-dependent Hamiltonian \cite{puzynin1999high,van2007accurate}.

The simplest approximation is the well-known Crank-Nicholson scheme. Many authors have noticed that the Crank-Nicholson scheme has the form of the canonical Cayley transform \cite{iserles2001cayley,havu2007cayley}. The Crank-Nicholson scheme and the schemes based on the diagonal Pad\'e approximation preserve the unitarity of the transition matrix, therefore they are used for systems with quadratic integrals \cite{diele1998cayley,van2007accurate}.
A general approach to constructing conservative one-step difference schemes can be considered. In this approach, two objects are subject to the definition: the generalized Cayley transform, which is given by an appropriate polynomial~$F(z)$, and the matrix~$Z$, which replaces the complex variable~$z$ in the generalized Cayley transform. Thus, schemes based on the Pad\'e approximation and the canonical Cayley transform are embedded in this approach, since in these cases the polynomial~$F(z)$ is specified a priori. In the most general setting, we can assume that the polynomial~$F(z)$ has complex coefficients. However, calculations showed that such a generalization does not allow decreasing the degrees of the polynomials in the spectral parameter. Therefore, we limited ourselves to polynomials~$F$ with real coefficients.

\subsection{Exponential integrators}

Let's introduce a fundamental solution $U(t,t_0)$ of the system
\begin{equation}
 \frac{d\,U(t,t_0)}{dt}=Q(t)\,U(t,t_0),\quad U(t_0,t_0)=I,
\end{equation}
where $I$ is a unit matrix. If the matrix $ Q $ does not depend on time, then the fundamental solution is the exponential $ U (t, t_0) = \exp\left((t-t_0) Q\right) $. Therefore, if $ Q $ depends on time, we can assume that the fundamental solution also has an exponential form. Indeed, in Magnus's work \cite {magnus1954exponential} the asymptotic representation of the fundamental solution in exponential form was found:
\begin{equation}\label{e^Omega}
 U(t,0)=e^{\Omega(t)},\quad \Omega(t)=\sum\limits_{k=0}^\infty\,\Omega_k(t),
\end{equation}
where the first terms of the Magnus expansion have the form
\begin{equation}
 \Omega_1(t)=\int\limits_0^tQ(t_1)dt_1,\quad
 \Omega_2(t)=\frac{1}{2}\int\limits_0^t dt_1 \int\limits_0^{t_1} dt_2 \left[Q(t_1),Q(t_2)\right],
\end{equation}
\begin{equation}\nonumber
 \Omega_3(t)=\frac{1}{6}\int\limits_0^t dt_1 \int\limits_0^{t_1} dt_2 \int\limits_0^{t_2} dt_3 \left(\left[Q(t_1),[Q(t_2),Q(t_3)]\right]
 +\left[Q(t_3),[Q(t_2),Q(t_1)]\right]\right).
\end{equation}

If we know how to calculate $Q$ at points inside the interval at each step of size $t$, then it is reasonable to approximate the integrals by suitable cubature formulas. In a situation where the matrix $Q$ is given on a uniform grid, various interpolation methods can be used.
Numerical experiments have shown the insufficiency of interpolation based on cubic splines and the efficiency of interpolation based on the Fourier transform \cite{chimmalgi2019fast}.
However, this interpolation is essentially nonlocal. Our numerical experiments with interpolation of smooth analytical signals using the Fourier transform showed that such interpolation gives values that coincide with the analytical ones with the accuracy of machine representation of numbers. Additional research is required to understand how this interpolation works for realistic signals. Therefore, another method was proposed for constructing exponential integrators using the Magnus formula.

If the matrix $Q(t)$ can be represented as a Taylor series with respect to a small parameter $\tau$
\begin{equation}
 Q(t+\tau)=\sum\limits_{k=0}^\infty\,\frac{\tau^k}{k!}Q^{(k)}(t),\quad Q^{(k)}(t)=\frac{d^kQ(t)}{dt^k},
\end{equation}
then substitution of this series into the Magnus formula with integration from $t-\tau/2$ to $t+\tau/2$ gives an approximation of the fundamental solution $U(t+\tau/2,t-\tau/2)$ with the required order of accuracy in the small parameter $\tau$.
We restrict ourselves to considering the approximation $E(t+\tau/2,t-\tau/2)$ up to the sixth order
in $\tau$
\begin{equation}\label{exp6}
 U(t+\tau/2,t-\tau/2)=E(t+\tau/2,t-\tau/2)+O\left(\tau^7\right),
\end{equation}
which has the form
\begin{equation}\label{e^Z}
 E(t+\tau/2,t-\tau/2)=e^{Z(t)},\quad Z(t)=\tau Z_1(t)+\tau^3Z_3(t)+\tau^5Z_5(t),
\end{equation}
where $Z_2=Z_6=0$ and the nonzero terms are
\begin{equation}\label{Z1eZ3e}
 Z_1=Q,\quad Z_3=\frac{1}{24}Q^{(2)}+\frac{1}{12}\left[Q^{(1)},Q\right],
\end{equation}
\begin{equation}\label{Z5e}
 Z_5=\frac{1}{1920}Q^{(4)}+\frac{1}{480}\left[Q^{(3)},Q\right]+\frac{1}{480}\left[Q^{(1)},Q^{(2)}\right]
+\frac{1}{720}\left[\left[Q^{(2)},Q\right],Q\right]+
\end{equation}
$$
 +\frac{1}{240}\left[\left[Q,Q^{(1)}\right],Q^{(1)}\right]
 +\frac{1}{720}\left[Q^3,Q^{(1)}\right]+\frac{1}{240}\left[QQ^{(1)}Q,Q\right].
$$

The expression $Z(t)$ contains derivatives from the first to the fourth order. To obtain a consistent finite-difference approximation, we express these derivatives on a uniform grid with a five-point stencil. We need the following central difference approximations of the derivatives, giving the maximum order of accuracy~\cite{fornberg1988generation}.

In the term $Z_5$, it is sufficient to approximate the derivatives with the 2nd order:
\begin{equation}
 Q^{(4)}(t)=\frac{Q_2-4Q_1+6Q_0-4Q_{-1}+Q_{-2}}{\tau^4}+O(\tau^2),
\end{equation}
\begin{equation}
 Q^{(3)}(t)=\frac{Q_2-2Q_1+2Q_{-1}-Q_{-2}}{2\tau^3}+O(\tau^2),
\end{equation}
\begin{equation}
 Q^{(2)}(t)=\frac{Q_1-2Q_0+Q_{-1}}{\tau^2}+O(\tau^2),
\end{equation}
\begin{equation}
 Q^{(1)}(t)=\frac{Q_1-Q_{-1}}{2\tau}+O(\tau^2),\quad Q_n=Q(t+n\tau).
\end{equation}

For $Q^{(2)}$ and $Q^{(1)}$ in the term $Z_3$ it is sufficient to use the 4th order approximation 
with a five-point stencil:
\begin{equation}
 Q^{(2)}(t)=\frac{-Q_2+16Q_1-30Q_0+16Q_{-1}-Q_{-2}}{12\tau^2}+O(\tau^4),
\end{equation}
\begin{equation}
 Q^{(1)}(t)=\frac{-Q_2+8Q_1-8Q_{-1}+Q_{-2}}{12\tau}+O(\tau^4).
\end{equation}
These finite-difference approximations can also be used for $Q^{(2)}$ and $Q^{(1)}$ in the term $Z_5$, which is equivalent to using the Lagrange interpolation polynomial of 4th degree in $\tau$ to approximate $Q(t+\tau)$ with a five-point stencil.

Formulas for the expansion of $\Omega(t)$ up to the 8th order in $\tau$ are given in~\cite{blanes2009magnus}. In this case, 18 nested commutators are added, which additionally contain $Q^{(6)}$ and $Q^{(5)}$. To obtain a consistent 8th order finite-difference scheme, at least a 7-point stencil must be used. 

If the matrix $Q$ is skew-Hermitian $Q^\dagger=-Q$, then the matrix $U(t,t_0)$ is unitary $U^{-1}(t,t_0)\equiv U(t_0,t)=U^\dagger(t,t_0)$. 
Obviously, all approximations of $Z(t)$ will also be skew-Hermitian and the finite-difference scheme based on the expansion will preserve the quadratic integral. In particular, if $Q(t)=-iH(t)$, where $H(t)$ is a Hermitian matrix, then we get the Schr\"odinger equation with the Hamiltonian $H(t)$ depending on the time $t$.

\subsection{Formulas for matrices of the second order}\label{appendix1}

Calculation of the matrix exponential $\exp(Z)$ in the general case is a rather complicated computational problem~\cite{moler2003nineteen}. However, for matrices of the 2nd and 3rd orders, the calculation of the matrix exponentials can be done analytically.
In this subsection, we consider the case of second-order matrices.

It is convenient to expand complex matrices of the second order in terms of the Pauli matrices
\begin{equation}\label{paulimatr}
 \sigma_0=\begin{bmatrix}1&0\\0&1\end{bmatrix}\equiv I,\quad \sigma_1=\begin{bmatrix}0&1\\1&0\end{bmatrix},\quad
 \sigma_2=\begin{bmatrix}0&-i\\i&0\end{bmatrix},\quad \sigma_3=\begin{bmatrix}1&0\\0&-1\end{bmatrix}.
\end{equation}

Let us consider a matrix $Z$, that has the following expansion in terms of the Pauli matrices $\sigma_k$:
\begin{equation}\label{Zpauli}
 Z=\begin{bmatrix}Z_{11}&Z_{12}\\Z_{21}&-Z_{11}\end{bmatrix}=\begin{bmatrix}z_3&z_1-iz_2\\z_1+iz_2&-z_3\end{bmatrix}=z_1\sigma_1+z_2\sigma_2+z_3\sigma_3,\quad z_k\in\mathbb{C}.
\end{equation}
This means that the matrix $Z$ is traceless: $\mbox{tr}(Z)=0$.
The characteristic equation of this matrix has the form
\begin{equation}
 \det(Z-\lambda \sigma_0)=\lambda^2-Z_{11}^2-Z_{12}Z_{21}=\lambda^2-z_1^2-z_2^2-z_3^2=0,
\end{equation}
from which we obtain an expression for the eigenvalues $\lambda_\pm=\pm\sqrt{z_1^2+z_2^2+z_3^2}$, and by the Hamilton-Cayley theorem it follows that $Z$ satisfies its characteristic equation
\begin{equation}
 Z^2=(z_1^2+z_2^2+z_3^2)\sigma_0
\end{equation}

The spectral decomposition of the matrix $Z=X\Lambda X^{-1}$ has the form
\begin{equation}\label{spectral_decomp}
 X=Z+\lambda\sigma_3,\quad \Lambda=\lambda\sigma_3,\quad X^{-1}=(2\lambda(z_3+\lambda))^{-1}X,
 \quad \lambda=\sqrt{z_1^2+z_2^2+z_3^2}.
\end{equation}

Let's take the analytic function $F(\lambda)$, which is given by its series.
Then the function $F$ of the matrix $Z$ takes the form
\begin{equation}
 F(Z)=F(X\lambda\sigma_3X^{-1})=X\left(F_c(\lambda)\sigma_0+F_s(\lambda)\sigma_3\right)X^{-1}
 =F_c(\lambda)\sigma_0+\frac{F_s(\lambda)}{\lambda}Z,
\end{equation}
where $F_c(\lambda)$ and $F_s(\lambda)$ are even and odd parts of the function $F(\lambda)$:
$$F_c(\lambda)=\frac{F(\lambda)+F(-\lambda)}{2},\quad F_s(\lambda)=\frac{F(\lambda)-F(-\lambda)}{2}.$$
In particular, for $F(\lambda)=e^\lambda$ we have
a compact formula for the exponential
\begin{equation}\label{e^Z=cs}
 e^Z=c(\lambda)\sigma_0+\frac{s(\lambda)}{\lambda}Z,
\end{equation}
where
\begin{equation}
 c(\lambda)=\cosh(\lambda),\quad s(\lambda)=\sinh(\lambda).
\end{equation}

For a rational function $F(\lambda)/G(\lambda)$ we obtain in a similar way the formula
\begin{equation}\label{FG}
\frac{F(Z)}{G(Z)}=\frac{F(\lambda)G(-\lambda)+F(-\lambda)G(\lambda)}{2G(\lambda)G(-\lambda)}\sigma_0+
\frac{1}{\lambda}\frac{F(\lambda)G(-\lambda)-F(-\lambda)G(\lambda)}{2G(\lambda)G(-\lambda)}Z.
\end{equation}

The generalized Cayley transform has the form $F(z)/F(-z)$.
Therefore, substituting $G(\lambda)=F(-\lambda)$ into the formula (\ref{FG}), we obtain the compact formula
for the generalized Cayley transform of a second-order matrix~$Z$ from (\ref{Zpauli})
\begin{equation}\label{PC}
\frac{F(Z)}{F(-Z)}=\frac{1}{F(\lambda)F(-\lambda)}\left(F_c(\lambda)\sigma_0+\frac{F_s(\lambda)}{\lambda}Z\right)^2
=c(\lambda)\sigma_0+\frac{s(\lambda)}{\lambda}Z,
\end{equation}
where $\lambda$ is the eigenvalue (\ref{spectral_decomp}) and 
the coefficients are
\begin{equation}\label{cs}
c(\lambda)=\frac{F^2(\lambda)+F^2(-\lambda)}{2F(\lambda)F(-\lambda)},\quad
s(\lambda)=\frac{F^2(\lambda)-F^2(-\lambda)}{2F(\lambda)F(-\lambda)}
\end{equation}
and satisfy the identity
\begin{equation}
c^2(\lambda)=1+s^2(\lambda).
\end{equation}

Obviously, $c(\lambda)$ and $s(\lambda)/\lambda$ are even functions of $\lambda$ and can be written as functions of $\lambda^2$. Also, if $F(\lambda)$ is a polynomial, then the degrees of the numerator and denominator for $c(\lambda)$ are the same, and the degree of the numerator of $s(\lambda)$ is less than the degree of the numerator $c(\lambda)$.

\subsection{Diagonal Pad\'e approximants}

For problems of the scattering theory, it is typical that the matrix $Q(t,\zeta)$ also depends on the spectral parameter $\zeta$ and it is necessary to find a solution of the equation (\ref{Psit}) for a large number of values of the parameters $\zeta$ to determine scattering data.
Therefore it was proposed \cite{wahls2013introducing,wahls2015fast} to represent the transition matrix $T_n=U(t_n+\tau/2,t_n-\tau/2)$ at each step $t_n$ as a rational function
\begin{equation}\label{Sd}
T_n=\frac{S_n(w)}{d_n(w)},
\end{equation}
where the matrix $S_n$ and the function $d_n$ are polynomials of the parameterization $w=w(\zeta)$. 
As a parameterization one can choose a function that transforms the space of the spectral parameter~$\zeta$ into a unit disc~$|w|\leq 1$.

In particular, if the spectral parameters~$\zeta$ lie in the upper half-plane, as in the ZS problem,
a linear fractional transformation that maps $\zeta=\xi+i\eta$, $\eta\geq 0$ into the unit disc $|w|\leq 1$ such that the point $i\beta$, $\beta>0$ goes to the point $w=\alpha$, $-1<\alpha<1$ inside the disc and the point $\zeta = 0$ goes to the point $w = 1$, has the form
\begin{equation}\label{w(z)}
 w(\zeta)=\frac{ih-\zeta}{ih+\zeta},\quad h=\frac{1+\alpha}{1-\alpha}\beta,\quad w(i\beta)=\alpha,\quad w(0)=1.
\end{equation}
The inverse transformation for (\ref{w(z)}) has the form
\begin{equation}\label{z(w)}
 \zeta(w)=ih\,\frac{1-w}{1+w}.
\end{equation}
This inverse transformation maps the points of the unit circle $w=\exp(i\theta)$, $\theta\in\mathbb{R}$ to points of the real axis $\zeta=\xi$, $\xi\in\mathbb{R}$ by the formula
\begin{equation}
 \xi=h\,\tan\left(\frac{\theta}{2}\right).
\end{equation}

Then, following (\ref{Sd}), the fundamental solution on the computational domain $[-\tau/2,\tau M+\tau/2]$  will be approximated by the product
\begin{equation}
U(\tau M+\tau/2,-\tau/2)\approx T(w)=\prod\limits_{n=0}^M T_n(w).
\end{equation}
The matrix $T(w)$ is a matrix polynomial with respect to $w$. Its coefficients can be found using fast algorithms for multiplying polynomials~\cite{pan2012structured, mcnamee2013numerical}. To compute the matrix $T(w)$ for different values of the free variable $w$ one can also use fast algorithms based on nonequispaced fast Fourier transform (NFFT)~\cite{Keiner2009a}.

To obtain a rational transition matrix $R$, one can use the Pad\'e approximation for the matrix exponential
\begin{equation}
 e^Z=R(Z)+O(\tau^k),\quad R(z)=\frac{F(z)}{G(z)},
\end{equation}
where $Z$ is a matrix depending on $\tau$, $F(z)$ and $G(z)$ are polynomials,
and the order of approximation $k$ must be no less than the order with which $\exp(Z)$ approximates the fundamental solution $U$. Using the adjugate matrix $\mbox{adj}(G(Z))$, the transition matrix will take the form (\ref{Sd}):
\begin{equation}\label{Sdadj}
 T_n=\frac{\mbox{adj}(G(Z_n))F(Z_n)}{\det(G(Z_n))}.
\end{equation}
For small dimensions 2 and 3, the inverse matrix $G^{-1}(Z_n)$ and/or the adjugate matrix $\mbox{adj}(G(Z_n))$ can be calculated analytically.

For equations with constant coefficients, rational approximations have been discussed for a long time. Moreover, the form of the polynomials must be consistent with the spectrum of the constant matrix $Q$. It is especially important for the stiff systems~\cite{blue1970rational}. For the Schr\"odinger equation with the time-dependent Hamiltonian~$H(t)$, difference schemes were constructed based on the diagonal Pad\'e approximation of the exponential. This ensures that the transition matrix is unitary. A general form of the diagonal Pad\'e approximation of the exponential is given in~\cite{baker1996pade} 
\begin{equation}\label{Paden}
 e^z=E_{n}(z)+O(z^{2n+1}),\quad E_{n}(z)=\frac{_1F_1(-n,-2n,z)}{_1F_1(-n,-2n,-z)},
\end{equation}
where $_1F_1(-n,-2n,z)$ is a confluent hypergeometric function that is reduced to a polynomial of degree $n$.
The first 4 diagonal Pade approximations have the form
\begin{equation}\label{Pade1}
E_1(z)=\frac{1+\frac{1}{2}z}{1-\frac{1}{2}z},
\end{equation}
\begin{equation}\label{Pade2}
E_2(z)=\frac{1+\frac{1}{2}z+\frac{1}{12}z^2}{1-\frac{1}{2}z+\frac{1}{12}z^2},
\end{equation}
\begin{equation}\label{Pade3}
E_3(z)=\frac{1+\frac{1}{2}z+\frac{1}{10}z^2+\frac{1}{120}z^3}{1-\frac{1}{2}z+\frac{1}{10}z^2-\frac{1}{120}z^3},
\end{equation}
\begin{equation}\label{Pade4}
E_4(z)=\frac{1+\frac{1}{2}z+\frac{3}{28}z^2+\frac{1}{84}z^3+\frac{1}{1680}z^4}
{1-\frac{1}{2}z+\frac{3}{28}z^2-\frac{1}{84}z^3+\frac{1}{1680}z^4}.
\end{equation}

For the 6th order exponential scheme (\ref{exp6}), one need to use the diagonal Pad\'e approximation, starting from the 3rd degree:
\begin{equation}\label{pade33}
 e^z=E_{3}(z)+O(z^7),
\end{equation}
where $E_3(z)$ is given in (\ref{Pade3}),
or more accurate approximations $E_{n}$ for $n\geq 3$ can be used.

For matrices $Q$ of the larger size, the polynomials $F_n(z)$ and $F_n(-z)$ 
can be factorized to represent a one-step difference scheme as a multi-step implicit scheme \cite{puzynin1999high,van2007accurate}.

\subsection{Generalized Cayley transform}

The Schr\"odinger equation and the ZS system for the real spectral parameter $\zeta=\xi\in\mathbb{R}$ preserve the quadratic integral. Therefore, we will construct transition matrices that also preserve this invariant.

For the approximate transition matrix $T$ to be unitary, it is sufficient that it has the form of the generalized Cayley transform
\begin{equation}\label{TF}
 T(Z)=\frac{F(Z)}{\overline{F}(-Z)},\quad F(z)=\sum\limits_{n=0}^\infty f_nz^n,\quad \overline{F}(z)=\sum\limits_{n=0}^\infty \bar{f}_nz^n
\end{equation}
where $F(z)$ is an analytic function of the complex argument $z$ such that
\begin{equation}
 F(z)\not\equiv \overline{F}(-z).
\end{equation}
For any real $y$ the generalized Cayley transform $w=F(iy)/\overline{F}(-iy)$ converts the imaginary axis $z=iy$ to the unit circle because $|w|=1$.
Further, we will consider only functions $F$ with real coefficients.
Obviously, the exponential $\exp(x)$ is an example of the generalized Cayley transform for $F(x)=\exp(x/2)$.

We will search the transition matrix $T$ in the form of the generalized Cayley transform up to 6th order in $\tau$. To do this, it suffices to consider the generalized Cayley transform in the form of a sixth-degree polynomial
\begin{equation}\label{F(Z)}
F(Z)=a_0I+a_1Z+a_2Z^2+a_3Z^3+a_4Z^4+a_5Z^5+a_6Z^6,
\end{equation}
and the expansion for $Z$
\begin{equation}\label{genZ}
Z=\tau Z_1+\tau^2Z_2+\tau^3Z_3+\tau^4Z_4+\tau^5Z_5+\tau^6Z_6,
\end{equation}
which starts with a first-order term in $\tau$.
The expansion of the fundamental solution $U(t+\tau/2,t-\tau/2)$ in $\tau$ has the form $I+\tau Q(t)$ in the main order, therefore $a_0$ and $a_1$ are not equal to zero. Without loss of generality, we can assume that $a_0=1$. Normalizing 
$Z$, we can set $a_1=1/2$:
\begin{equation}\label{cayley6}
T(Z)=\frac{I+\frac{1}{2}Z+a_2Z^2+a_3Z^3+a_4Z^4+a_5Z^5+a_6Z^6}{I-\frac{1}{2}Z+a_2Z^2-a_3Z^3+a_4Z^4-a_5Z^5+a_6Z^6}.
\end{equation}
For $a_2=a_3=a_4=a_5=a_6$ we obtain exactly the canonical Cayley transform 
\begin{equation}\label{canCayley}
E_1(Z)=\frac{I+\frac{1}{2}Z}{I-\frac{1}{2}Z}.   
\end{equation}

Formulas for $Z$ from (\ref{genZ}) in the 6th order general scheme (\ref{cayley6}) take the form: $Z_2=Z_4=Z_6=0$,
\begin{equation}\label{Z1}
Z_1=Q,
\end{equation}
\begin{equation}\label{Z3}
Z_3=\frac{1}{24}Q^{(2)}+\frac{1}{12}\left[Q^{(1)},Q\right]+k_1Q^3,
\end{equation}
\begin{equation}\label{Z5}
Z_5=\frac{1}{1920}Q^{(4)}+\frac{1}{480}\left[Q^{(3)},Q\right]+\frac{1}{480}\left[Q^{(1)},Q^{(2)}\right]
+\frac{1}{240}\left[\left[Q,Q^{(1)}\right],Q^{(1)}\right]
\end{equation}
$$
+k_2\left[\left[Q^{(2)},Q\right],Q\right]+k_3QQ^{(2)}Q
+k_4\left[Q^3,Q^{(1)}\right]+\frac{1}{240}\left[QQ^{(1)}Q,Q\right]+k_5Q^5,
$$
where
\begin{equation}
 k_1=a_2-2a_3-\frac{1}{12},\quad k_2=\frac{1}{24}\left(k_1+\frac{1}{30}\right),\quad k_3=\frac{k_1}{8},
 \quad k_4=-\frac{1}{12}\left(k_1-\frac{1}{60}\right)
\end{equation}
\begin{equation}
 k_5=\frac{1}{120}-\frac{a_2}{4}+\frac{a_3}{2}+2a_2^2-10a_2a_3+12a_3^2+a_4-2a_5.
\end{equation}
Another notation of $k_5$ through $k_1$ has the form
\begin{equation}
 k_5=2k_1\left(k_1-a_3+\frac{1}{24}\right)-\frac{1}{6}\left(a_3-\frac{1}{120}\right) +a_4-2a_5.
\end{equation}

These formulas do not contain the coefficient $a_6$, since, as in the exponential expansion, 
the matrix $Z$ has only the odd powers of $\tau$.
There are four arbitrary coefficients: $a_2$, $a_3$, $a_4$, $a_6$. Moreover $a_4$ and $a_5$ are included only in the coefficient at $Q^5$ in the form of a linear combination $a_4-2a_5$.

For the third order diagonal Pad\'e approximation (\ref{pade33}), 
i.e. for $a_2=1/10$, $a_3=1/120$, $a_4=a_5=0$,
the matrix $Z$ coincides with the matrix (\ref{Z1eZ3e})-(\ref{Z5e}) for the exponential scheme (\ref{e^Z}). Thus, the general schemes (\ref{cayley6}) contain the third order Pad\'e approximation for the 6th order exponential scheme (\ref{e^Z}).

Arbitrariness in the choice of coefficients $a_k$, $k=2,3,4,5$, can be used in several ways.

First, to zero out the maximum number of terms in $Z$, we have to set $k_1 = 0$, then two terms  $k_1$ and $k_3$ are canceled.
Putting $a_4=a_5=0$, to decrease the degree of the polynomial $F(z)$, we will zero out the coefficient $k_5$ for $a_3=1/120$. As a result, we get $a_2=1/10$. Therefore, this case coincides with the 3rd order Pad\'e approximation (\ref{Pade3}).

Second, to obtain the minimum degree of a polynomial, we put $a_2=a_3=a_4=a_5=0$. Then we get the canonical Cayley transform (\ref{canCayley}), and the matrix $Z$ will be determined by the coefficients
\begin{equation}\label{kk}
 k_1=-\frac{1}{12},\quad k_2=-\frac{1}{480},\quad k_3=-\frac{1}{96},\quad k_4=k_5=\frac{1}{120}.
\end{equation}

Let us consider the question: how can we choose the coefficients $a_k$ so that the polynomials $F(z)$ and $F(-z)$ have a common root that can be canceled in a rational expression (\ref{cayley6})? For two polynomials $F(z)$ and $F(-z)$ have a common root, it is enough that their resultant is equal to zero. 
Calculations for polynomials of the 5th degree show that, under the condition
\begin{equation}
a_2=2a_3,\quad a_4=2a_5,
\end{equation}
the maximum reduction occurs up to polynomials of the first degree, i.e. to the canonical Cayley transform~(\ref{canCayley}), and the matrix~$Z$ is determined by the coefficients~(\ref{kk}). 
Another case of reduction is to a polynomial of the 3rd degree, but this does not zero the coefficient $k_5$, so a polynomial of a higher degree is obtained than for the canonical Cayley transform~(\ref{canCayley}).

\subsection{Conditions of applicability for schemes}

The transition matrix $T$ is close to the unit matrix $I$ for sufficiently small $\tau$. The approximate transition matrix in exponential form (\ref{e^Z}) satisfies this property for any $\tau$. If the matrix $Z$ has a simple structure and $\lambda_k$ is a set of eigenvalues, then the polynomials $F(\lambda_k)$ and $F(-\lambda_k)$ from (\ref{TF}) have to be far from their zeros. For Pad\'e approximation of an exponential function, zeros and poles are well studied \cite{baker1996pade,saff1978zeros}. 
For small orders of the generalized Pad\'e transform (\ref{cayley6}), the zeros of the numerator and denominator can be found numerically.
If the root with the minimum modulus of the polynomial $F(z)$ is equal to $z_*$, then the condition of applicability of the difference scheme can be written in the form
\begin{equation}
|z_*|>|\lambda_k|
\end{equation}
for all $k$.

\section{Zakharov-Shabat system}\label{ZS}

In this section, we will consider a modified ZS system with the matrix
\begin{equation}\label{Qqr}
 Q(t,\zeta)=\begin{bmatrix}-i\zeta&q(t)\\r(t)&i\zeta\end{bmatrix}.
\end{equation}
For different functions $q(t)$ and $r(t)$, the modified ZS system corresponds to the direct spectral problem for some nonlinear equations. A list of such equations is given in \cite{ablowitz1981solitons,wahls2013introducing,wahls2015fast}. In addition, the ZS system is used to describe the integrable generalizations of the NLSE, 
which can be used to describe the pulse propagation in optical fibers \cite{hasegawa1995solitons}. Using the general theory from the previous section, we construct three sixth-order difference schemes for this system. The schemes based on the diagonal Pad\'e approximation and the Cayley transform allow the use of fast algorithms to solve the direct spectral problem for a large number of values of the spectral parameter $\zeta$.

\subsection{Demo example}

We will consider the Crank-Nicholson scheme for the system (\ref{Psit}) with the matrix (\ref{Qqr}) to demonstrate the use of the fast algorithm. This scheme, like several other schemes, was considered in \cite{wahls2013introducing,wahls2015fast}. We chose it because the Crank-Nicholson scheme is a prototype for schemes based on the diagonal Pad\'e approximation and the generalized Cayley transform, and the formulas for it have the most compact form.

The transition matrix $T$ for the exponential scheme of the second order of accuracy has the form
\begin{equation}\label{bo}
 T=e^{Z}=\cosh(\lambda)\sigma_0+\frac{\sinh(\lambda)}{\lambda}Z,\quad Z=\tau Q=\begin{bmatrix}0&\tau q\\\tau r&0\end{bmatrix}-i\tau\zeta\begin{bmatrix}1&0\\0&-1\end{bmatrix},\quad \lambda=\tau\sqrt{qr-\zeta^2}.
\end{equation}
The matrix $Q$ has the inverse time dimension, therefore the matrix $Z$ is dimensionless and it is necessary to use dimensionless combinations $\tilde{q}=\tau q$, $\tilde{r}=\tau r$, and $z=\tau\zeta$. The scheme (\ref{bo}) was proposed for the ZS system in~\cite{boffetta1992computation}.

The first order diagonal Pad\'e approximation $E_1(z)$ approximates the exponent $\exp(z)$ with the second order of accuracy and has the form (\ref{Pade1}), therefore the corresponding transition matrix $T$ is written as
\begin{equation}\label{p2z}
 T(z)=
 \frac{S(z)}{d(z)},\quad S(z)=\left(1+\frac{1}{4}\tilde{q}\tilde{r}-\frac{1}{4}z^2\right)
 \begin{bmatrix}1&0\\0&1\end{bmatrix}+
 \begin{bmatrix}0&\tilde{q}\\\tilde{r}&0\end{bmatrix}
 -iz\begin{bmatrix}1&0\\0&-1\end{bmatrix},\quad d(z)=1-\frac{1}{4}\tilde{q}\tilde{r}+\frac{1}{4}z^2.
\end{equation}
Further within Section~\ref{ZS}, to simplify the notation, we will remove the wave over $q$ and $r$.

Let us perform a one-to-one conformal linear fractional transformation of the unit disc $|w|\leq 1$ into the upper half-plane $\mbox{Re}\,z\geq 0$
\begin{equation}\label{z(w)=ih}
 z(w)=ih\,\frac{1-w}{1+w},
\end{equation}
where $h>0$ is a parameter (\ref{w(z)}). 
Substituting (\ref{z(w)=ih}) into (\ref{p2z}) we get expressions for $S$ and $d$ in terms of $w$
\begin{equation}\label{p2w}
P_2(w)=
\frac{S(w)}{d(w)},\quad S(w)=S_0(w)\begin{bmatrix}1&0\\0&1\end{bmatrix}+
S_{12}(w)\begin{bmatrix}0&q\\r&0\end{bmatrix}
+S_3(w)\begin{bmatrix}1&0\\0&-1\end{bmatrix},
\end{equation}
where
\begin{equation}
d(w)=(4-h^2-qr)w^2+(8+2h^2-2qr)w+(4-h^2-qr),
\end{equation}
\begin{equation}
S_0(w)=(4+h^2+qr)w^2+(8-2h^2+2qr)w+(4+h^2+qr),
\end{equation}
\begin{equation}
S_{12}(w)=4(w+1)^2,\quad S_3(w)=4h(1-w^2).
\end{equation}
Now all calculations of matrix polynomials $S(w)$ and $d(w)$ will be performed for $|w|\leq 1$.

\subsection{Exponential and rational schemes of 6th order of accuracy}
Substituting the matrix $Q$ from (\ref{Qqr}) into general formulas (\ref{Z1eZ3e})-(\ref{Z5e}), we obtain the matrix $Z$ for the exponential scheme (\ref{e^Z}) in the compact form (\ref{e^Z=cs}):
\begin{equation}\label{Z11exp}
 Z_{11}=-Z_{22}=\frac{1}{180}\left(rq^{(1)}-qr^{(1)}\right)z^2-
 i\left(1-\frac {rq^{(2)}+qr^{(2)}}{360}+\frac{q^{(1)}r^{(1)}}{60}\right)z+
\end{equation}
$$+\frac{1}{180}\left(15-qr\right)\left(rq^{(1)}-qr^{(1)}\right)+
\frac{1}{480}\left(rq^{(3)}-qr^{(3)}+q^{(1)}r^{(2)}-r^{(1)}q^{(2)}\right),$$
\begin{equation}
 Z_{12}=\frac{iq^{(1)}}{90}z^3-\frac{q^{(2)}}{180}z^2+
 i\left(\frac{q^{(1)}}{6}+\frac{q^{(3)}}{240}-\frac{qrq^{(1)}}{90}\right) z+
\end{equation}\label{Z12exp}
$$+q+\frac{q\left(rq^{(2)}-qr^{(2)}\right)}{360}+\frac{\left(qr^{(1)}-rq^{(1)}\right)q^{(1)}}{120}
+\frac{q^{(2)}}{24}+\frac{q^{(4)}}{1920}, $$
\begin{equation}\label{Z21exp}
 Z_{21}=-\frac{ir^{(1)}}{90}z^3-\frac{r^{(2)}}{180}z^2-
 i\left(\frac{r^{(1)}}{6}+\frac{r^{(3)}}{240}-\frac{qrr^{(1)}}{90}\right) z+
\end{equation}
$$+r+\frac{r\left(qr^{(2)}-rq^{(2)}\right)}{360}+\frac{\left(rq^{(1)}-qr^{(1)}\right)r^{(1)}}{120}
+\frac{r^{(2)}}{24}+\frac{r^{(4)}}{1920}. $$

Sixth-order rational approximations are constructed using this expression for $Z$ and general formulas for the Pad\'e approximation (\ref{Paden}). 
The transition matrix $T=E_n(Z)$ has the form of (\ref{PC}) 
\begin{equation}\label{TPadeCayley}
T=E_n(Z)=c_n(\lambda)\sigma_0+\frac{s_n(\lambda)}{\lambda}Z  
\end{equation}
and for $n=3,4$ the coefficients $c_n(\lambda)$ and $s_n(\lambda)$ are
\begin{equation}\label{c3s3}
c_3(\lambda)=\frac{1+\frac{9}{20}\lambda^2+\frac{11}{600}\lambda^4+\frac{1}{14400}\lambda^6}
{1-\frac{1}{20}\lambda^2+\frac{1}{600}\lambda^4-\frac{1}{14400}\lambda^6},\quad
\frac{s_3(\lambda)}{\lambda}=\frac{1+\frac{7}{60}\lambda^2+\frac{1}{600}\lambda^4}
{1-\frac{1}{20}\lambda^2+\frac{1}{600}\lambda^4-\frac{1}{14400}\lambda^6},
\end{equation}
\begin{equation}\label{c4s4}
c_4(\lambda)=\frac{1+\frac{13}{28}\lambda^2+\frac{289}{11760}\lambda^4+\frac{19}{70560}\lambda^6+\frac{1}{2822400}\lambda^8}
{1-\frac{1}{28}\lambda^2+\frac{3}{3920}\lambda^4-\frac{1}{70560}\lambda^6+\frac{1}{2822400}\lambda^8},\quad
\frac{s_4(\lambda)}{\lambda}=\frac{1+\frac{11}{84}\lambda^2+\frac{37}{11760}\lambda^4+\frac{1}{70560}\lambda^6}
{1-\frac{1}{28}\lambda^2+\frac{3}{3920}\lambda^4-\frac{1}{70560}\lambda^6+\frac{1}{2822400}\lambda^8}.
\end{equation}

\subsection{Scheme for the canonical Cayley transform}
The transition matrix $T$ has the form of (\ref{TPadeCayley}) for $n=1$.
The matrix $Z$ of the 6th order scheme for the canonical Cayley transform  (\ref{canCayley}) has the form
\begin{equation}\label{Z11cayley}
Z_{11}=-Z_{22}=-\frac {i}{120}z^5+\frac {i}{60}(qr-5)z^3-\frac {1}{80} \left(qr^{(1)}-rq^{(1)} \right)z^2
\end{equation}
$$
-\frac {i}{480} \left( 4\,q^2r^2-40\,qr-3\,qr^{(2)}+8\,q^{(1)}r^{(1)}-3\,rq^{(2)}+480 \right) z$$
$$
+{\frac {1}{480}}
\left(6\,q^2rr^{(1)}-6\,qq^{(1)}r^2-40\,qr^{(1)}-r^{(3)}q+40\,q^{(1)}r+q^{(1)}r^{(2)}-r^{(1)}q^{(2)}+q^{(3)}r\right),$$
\begin{equation}\label{Z12cayley}
Z_{12}=\frac{q}{120}\,z^4+\frac{i}{40}q^{(1)}\,z^3+ \left(\frac{q}{12}-\frac{q^{(2)}}{480}-\frac{q^2r}{60}\right)z^2+i\left( \frac{1}{6}q^{(1)}-\frac{1}{40}qrq^{(1)}+\frac{1}{240}q^{(3)}\right)z
\end{equation}
$$
-\frac {(q^{(1)})^2r}{120}+\frac{q^{(1)}r^{(1)}q}{120}-\frac{qrq^{(2)}}{240}-\frac{q^2r}{12}+\frac{q^3r^2}{120}+\frac{q^{(4)}}
{1920}+\frac{q^{(2)}}{24}+q-\frac{q^2r^{(2)}}{160},$$
\begin{equation}\label{Z21cayley}
Z_{21}=\frac{r}{120}\,z^4-\frac{i}{40}r^{(1)}\,z^3+\left(\frac{r}{12}-\frac{r^{(2)}}{480}-\frac{qr^2}{60}\right)z^2-i\left( \frac{1}{6}r^{(1)}-\frac{1}{40}qrr^{(1)}+\frac{1}{240}r^{(3)}\right)z
\end{equation}
$$
-\frac {(r^{(1)})^2q}{120}+\frac{q^{(1)}r^{(1)}r}{120}-\frac{qrr^{(2)}}{240}-\frac{r^2q}{12}+\frac{r^3q^2}{120}+\frac{r^{(4)}}
{1920}+\frac{r^{(2)}}{24}+r-\frac{r^2q^{(2)}}{160}.$$

The functions $c_1(\lambda)$ and $s_1(\lambda)$ are calculated by the formulas
\begin{equation}\label{c1s1}
c_1(\lambda)=\frac{1+\frac{1}{4}\lambda^2}{1-\frac{1}{4}\lambda^2},\quad \frac{s_1(\lambda)}{\lambda}=\frac{1}{1-\frac{1}{4}\lambda^2}.
\end{equation}

\section{Numerical experiments}\label{num}

Let us consider numerical experiments for the constructed schemes using the example of the direct spectral problem for the ZS system.
Let $q=q(t,z)$ is a slow-varying complex optical field envelope propagating along an ideally lossless and noiseless fiber. The evolution of the pulse $q$ is described by the
standard NLSE
\begin{equation}\label{nlse}
i\frac{\partial q}{\partial z}+\frac{\sigma}{2}\frac{\partial^2 q}{\partial t^2}+|q|^2q=0,
\end{equation}
where the variable $z$ is the distance along the optical fiber, $t$ is a time variable; $\sigma =-1$ and $\sigma =1$ corresponds to the normal and anomalous dispersion in the fiber, respectively~\cite{Agrawal}.

The Nonlinear Fourier Transform allows to transform any signal $q(t)$, which decays rapidly for $t\to \pm \infty$, into nonlinear Fourier spectrum. It is defined by the solution of the ZS problem
\begin{equation}\label{psit}
\frac{d \Psi(t)}{dt}=Q(t){\Psi}(t),\quad
{\Psi}(t)=\begin{bmatrix}\psi_1(t)\\\psi_2(t)\end{bmatrix},\quad
Q(t)=\begin{bmatrix}-i\zeta&q(t)\\-\sigma q^*(t)&i\zeta\end{bmatrix},
\end{equation}
where $\Psi(t)$ is a complex vector function of a real argument $t$, $\zeta \in \mathbb{C}$ is a spectral parameter, $q(t)=q(t,z_0)$ for any fixed $z_0$.

Under the assumption that $q(t)$ decays rapidly when $t\to \pm \infty $, the specific solutions (Jost functions) for ZS problem (\ref{psit}) can be derived as
\begin{equation}\label{psi0}
\Psi =
\begin{bmatrix}
\psi_{1}\\\psi_{2}
\end{bmatrix} =
\begin{bmatrix}
e^{-i\zeta t}\\0
\end{bmatrix}[1+o(1)],\quad t\to-\infty,
\end{equation}
and
\begin{equation}\label{psi0right}
\Phi =\begin{bmatrix}
\phi_{1}\\\phi_{2}
\end{bmatrix} = 
\begin{bmatrix}
0\\e^{i\zeta t}
\end{bmatrix}[1+o(1)],\quad t\to \infty,
\end{equation}
Then we obtain the Jost scattering coefficients $a(\xi )$ and $b(\xi )$ as follows:
\begin{equation}\label{ab}
a(\xi)=\lim_{t\to\infty}\,\psi_1(t,\xi)\,e^{i\xi t},\quad b(\xi)=\lim_{t\to\infty}\,\psi_2(t,\xi)\,e^{-i\xi t}.
\end{equation}
The functions $a(\xi )$ and $b(\xi )$ can be extended to the upper half-plane $\xi \to \zeta$, where $\zeta $ is a complex number with the positive imaginary part~\cite{ablowitz1981solitons}. The spectral data of ZS problem (\ref{psit}) are determined by $a(\zeta)$ and $b(\zeta)$ in the following way:\\
(1) the continuous spectrum is determined by the reflection coefficient
$r(\xi)=b(\xi)/a(\xi)$, $\xi\in\mathbb{R}$.\\
(2) in the case of $\sigma =1$, the discrete spectrum $\{\zeta_k\}$, $k=\overline{0,K-1}$ is defined by $K$ zeros of $a(\zeta)=0$, and corresponding phase coefficients are defined as
$$r(\zeta_k)=\left.\frac{b(\zeta)}{a'(\zeta)}\right|_{\zeta=\zeta_k},\quad\mbox{where}\quad a'(\zeta)=\frac{da(\zeta)}{d\zeta};$$

The ZS system (\ref{psit}) conserves the quadratic invariant 
$H=|\psi_1|^2+\sigma|\psi_2|^2$ for real spectral parameters $\zeta=\xi$. In particular,
\begin{equation} \label{QuadInvariant}
H(\xi) = |a(\xi)|^2 + \sigma |b(\xi)|^2 = 1.
\end{equation}
In addition, the continuous spectrum energy 
\begin{equation} \label{Econt}
E_{c} =-\frac{1}{\pi } \int _{-\infty }^{\infty } \ln |a(\xi )|^{2} d\xi
\end{equation}
also conserves. The details of the conservative properties of the ZS system can be found in \cite{medvedev2019exponential,medvedev2020exponential}.

Summing up, we solve a linear system of the form~(\ref{psit}) with the matrix~$Q(t)$ linearly dependent on the complex function~$q(t)$. The numerical implementation of the continuous function~$ q(t)$ is a discrete function~$q_n=q(t_n)$, which is defined at the integer nodes $t_n$ of the uniform grid with the step $\tau$. Since we are considering a finite time interval, we will solve the problem on the interval $[-L,L]$ with the total number of points equal to $M+1$. In this case, the grid step is $\tau=2L/M$ and $t_n=-L+\tau n$, where $n=0,...,M$.

We replace the original system (\ref{psit}) on each subinterval $(t_n-\tau/2,t_n+\tau/2)$ with an approximate system with constant coefficients
\begin{equation}\label{T0}
\Psi(t_n+\tau/2)=T_n\Psi(t_n-\tau/2),
\end{equation}
where $T_n$ is a transition matrix from the layer $n-\frac{1}{2}$ to the layer $n+\frac{1}{2}$.

The spectral data are finally defined by
\begin{equation}\label{ab_compute}
a(\zeta)= \psi_1(L-\tau/2,\zeta)\,e^{i\zeta (L-\tau/2)},\quad 
b(\zeta)= \psi_2(L-\tau/2,\zeta)\,e^{-i\zeta (L-\tau/2)}.
\end{equation}

To compute the transition matrix $T_n = \exp(Z)$ (\ref{e^Z=cs}),
we need to find Pauli coefficients $z_1$, $z_2$, $z_3$ (\ref{Zpauli}) of the decomposition of the matrix $Z$ using Pauli matrices (\ref{paulimatr}). The Pauli coefficients $z_1$, $z_2$, $z_3$ are polynomials of a variable $z = \tau \zeta$. To optimize the calculations we compute coefficients of these polynomials for each grid node~$t_n$ at the preliminary stage. This procedure allows us to speed up the calculations since we will not need to compute the same coefficients for each value of the spectral parameter~$\zeta$. For a large number of spectral parameters, it gives a significant advantage. At the next stage, the problem (\ref{T0}) is solved for each value of the spectral parameters.

Here we consider the 6th order exponential scheme ES6 with the transition matrix $T_n = \exp(Z)$ (\ref{e^Z=cs}), where the matrix $Z$ is defined by (\ref{Z11exp})--(\ref{Z21exp}). Hyperbolic sine and cosine for the exponential scheme are calculated directly. 

We compare ES6 with two 6th order schemes ES6\_Pade3 and ES6\_Pade4 based on the diagonal Pad\'e approximations of the 3rd and 4th order.
The transition matrix $T_n$ is defined by the general formula (\ref{TPadeCayley}) for $n=3$ and $4$. In particular, the formula (\ref{Pade3}) is used for ES6\_Pade3 and the formula (\ref{Pade4}) is for ES6\_Pade4. The matrix $Z$ is defined by the same formulas (\ref{Z11exp})--(\ref{Z21exp}). The corresponding coefficients $c(\lambda)$ and $s(\lambda)$ are computed by (\ref{c3s3}) for ES6\_Pade3 and (\ref{c4s4}) for ES6\_Pade4.

We also consider the 6th order scheme ES6\_Cayley based on the Cayley transform with the transition matrix~$T_n$ (\ref{TPadeCayley}) for $n=1$ (\ref{Pade1}), where the matrix~$Z$ is defined by (\ref{Z11cayley})--(\ref{Z21cayley}) and the approximation (\ref{c1s1}) is used to find the corresponding coefficients $c(\lambda)$ and $s(\lambda)$.

The fast variants of the schemes (FES6\_Pade3 and FES6\_Pade4) were implemented 
based on the FNFT software library~\cite{Wahls2018}. We used NFFT3 library~\cite{Keiner2009a} to compute continuous spectrum by these schemes.
The fast variant of the scheme ES6\_Cayley turned out to be very inaccurate, so we do not present it in the figures.
Optimal values of the parameter $h$ (\ref{w(z)}) were chosen empirically as $h=11$ for the FES6\_Pade3 and $h=15$ for the FES6\_Pade4. But we should note that there are exist maximal critical value for these parameters ($h=11.65$ and $h=15.57$, correspondingly) when the schemes still work. For such parameters, the schemes demonstrate an exponential decrease in error for large step sizes $\tau$ but accumulate a lot of computational errors for small step sizes.

If we omit the term with $Z_5$ in formula (\ref{e^Z}), then we obtain an exponential scheme of the 4th order.
The numerical results for the exponential 4th order scheme ES4 can be found in the recent papers~\cite{medvedev2019exponential,medvedev2020exponential,medvedev2020conservative}.

We compared the forementioned exponential 6th order schemes with the $\text{CF}^{[6]}_4$ scheme \cite{chimmalgi2019fast, blanes2017high}. This is a
commutator-free quasi-Magnus (CFQM) exponential integrator with complex coefficients. Because of complex coefficients this scheme does not conserve the quadratic invariant. For the same reason it can not be made fast using the splitting method \cite{Prins2018a_sample}. The $\text{CF}^{[6]}_4$ scheme requires interpolation in two additional nodes for each subinterval. A sufficient result is given by interpolation based on the Fourier transform~\cite{chimmalgi2019fast}. 

We have applied the interpolation procedure to our schemes to provide a correct comparison at the grids with the same number of nodes.
The $\text{CF}^{[6]}_4$ scheme has shown almost the same accuracy as the ES6 scheme, but the running time of the $\text{CF}^{[6]}_4$ is longer.
Despite the fact that the matrices, composed the transition matrix in the $\text{CF}^{[6]}_4$, are much simpler, than the ones in the ES6, the calculations of four matrix exponentials in the $\text{CF}^{[6]}_4$ requires more time, then the computing the one matrix exponential in the ES6.
 The numerical experiments have also confirmed that the $\text{CF}^{[6]}_4$ scheme does not conserve the quadratic invariant.
Here we only present the graphs of the schemes without interpolation.

In \cite{chimmalgi2019fast} the fast sixth order scheme $\text{FCF\_RE}^{[4]}_2$ is also presented.
This scheme was constructed by integrating Richardson extrapolation into the fast fourth-order scheme $\text{FCF}^{[4]}_2$. The initial $\text{CF}^{[4]}_2$ scheme is a CQFM exponential integrator consisting of two exponentials and requiring interpolation in two additional nodes per subinterval. The $\text{CF}^{[4]}_2$ scheme conserves the quadratic invariant, but its fast variant $\text{FCF}^{[4]}_2$ does not, as well as the $\text{FCF\_RE}^{[4]}_2$.
We do not consider the $\text{FCF\_RE}^{[4]}_2$ in the current paper since we
believe that Richardson extrapolation is an improvement that can be applied to other schemes as well. The application of Richardson extrapolation to the exponential schemes and their subsequent comparison, in particular with the
$\text{FCF\_RE}^{[4]}_2$, is undoubted of interest and will be done in our future works.

For numerical experiments we used a conventional model signal in the form of a chirped hyperbolic secant
$q(t) = A[\mbox{sech}(t)]^{1+iC}$
with the following parameters: $A = 5.2$, $C = 4$ for both anomalous and normal dispersion.
The detailed analytical expressions of the spectral data for this type of potentials can be found in~\cite{medvedev2019exponential,medvedev2020exponential}.

To find the numerical errors of calculating the continuous spectrum energy~$E_c$ (\ref{Econt}), the quadratic invariant $H(\xi)$ (\ref{QuadInvariant}), the phase coefficients $r(\zeta_k)$, and the scattering coefficients $a(\zeta_k)$, $b(\zeta_k)$ at the eigenvalues~$\zeta_k$ we use the formula
\begin{equation}\label{err}
\mbox{err}[\phi]\!=\!\frac{|\phi^{comp} - \phi^{exact}|}{\phi_0}, \quad
\phi_0=
\begin{cases}
\phi^{exact},\!\mbox{ if }\!|\phi^{exact}|\!>\!1\\
1, \mbox{otherwise},
\end{cases}
\end{equation}
where $\phi$ can represent $E_c$, $H(\xi)$, $r(\zeta_k)$, $a(\zeta_k)$ or $b(\zeta_k)$.

For the continuous spectrum we calculate the root mean squared error
\begin{equation}\label{MSE}
RMSE[\phi]=\sqrt{\frac{1}{N}\sum_{j=1}^{N}\frac{|\phi^{comp}(\xi_j) - \phi^{exact}(\xi_j)|^2}{|\phi_0(\xi_j)|^2}},
\end{equation}
\begin{equation}\nonumber
\phi_0 = 
\begin{cases}
\phi^{exact}(\xi_j), \mbox{ if } |\phi^{exact}(\xi_j)| > 1\\
1, \mbox{ otherwise},
\end{cases}
\end{equation}
where $\phi$ can represent $a(\xi)$, $b(\xi)$, $r(\xi)$ or $H(\xi)$. Here we assume the spectral parameter $\xi \in [-20, 20]$ with the total number of points $N$ that equal to the number of points $M$ of the signal discretization. 

All calculations were performed on a single core of the Intel\textsuperscript{\textregistered} Core\textsuperscript{TM} i5-9600K processor with a frequency of 4.6 GHz. All algorithms were implemented using C++ language and compiled by Intel\textsuperscript{\textregistered} C++ Compiler 19.1.

\begin{figure}[htpb]
	\centering
	\includegraphics[width=13cm]{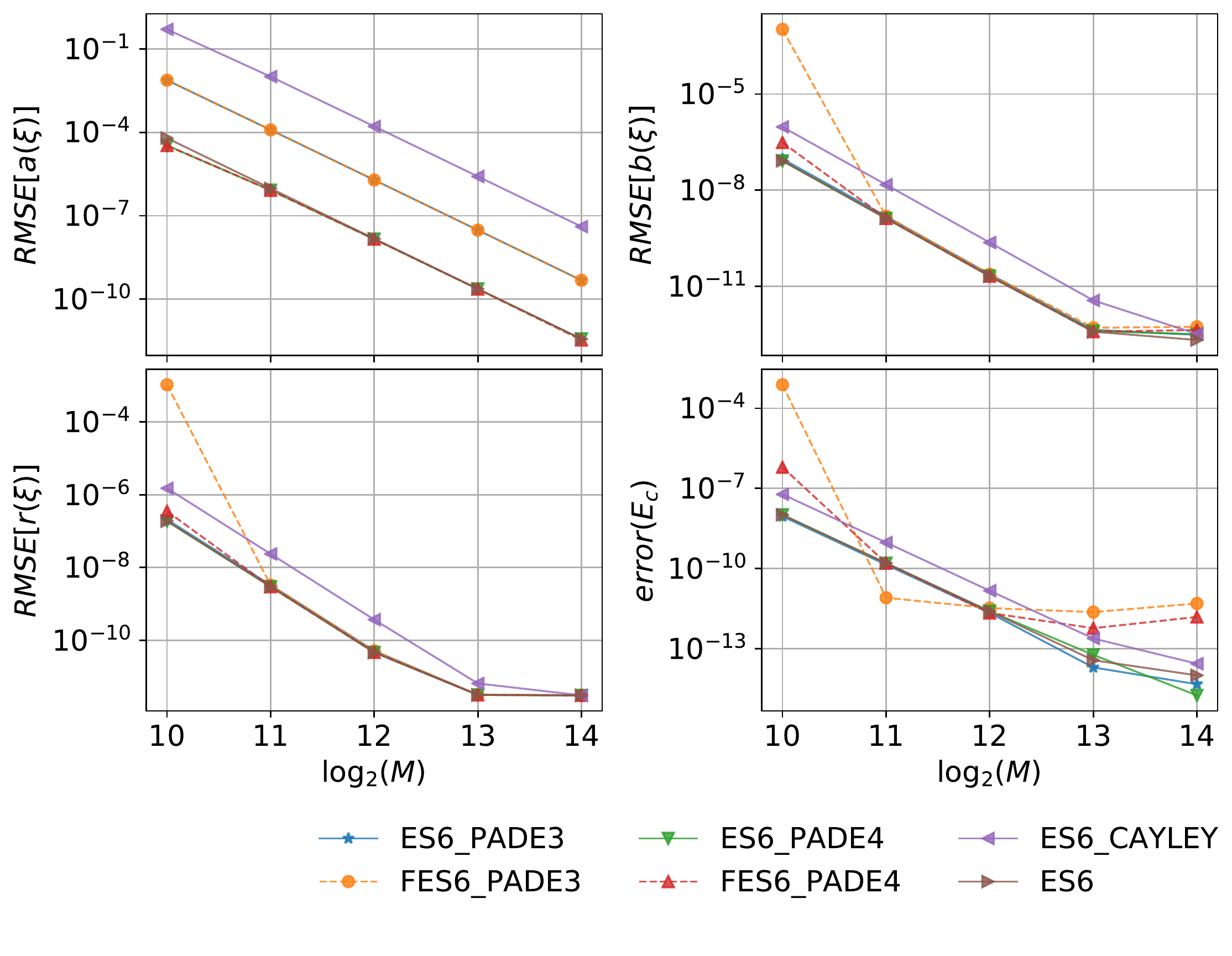}
	\caption{Continuous spectrum errors in the case of anomalous dispersion $\sigma=1$.	}
	\label{fig:error1}
\end{figure}
\begin{figure}[htpb]
	\centering
	\includegraphics[width=13cm]{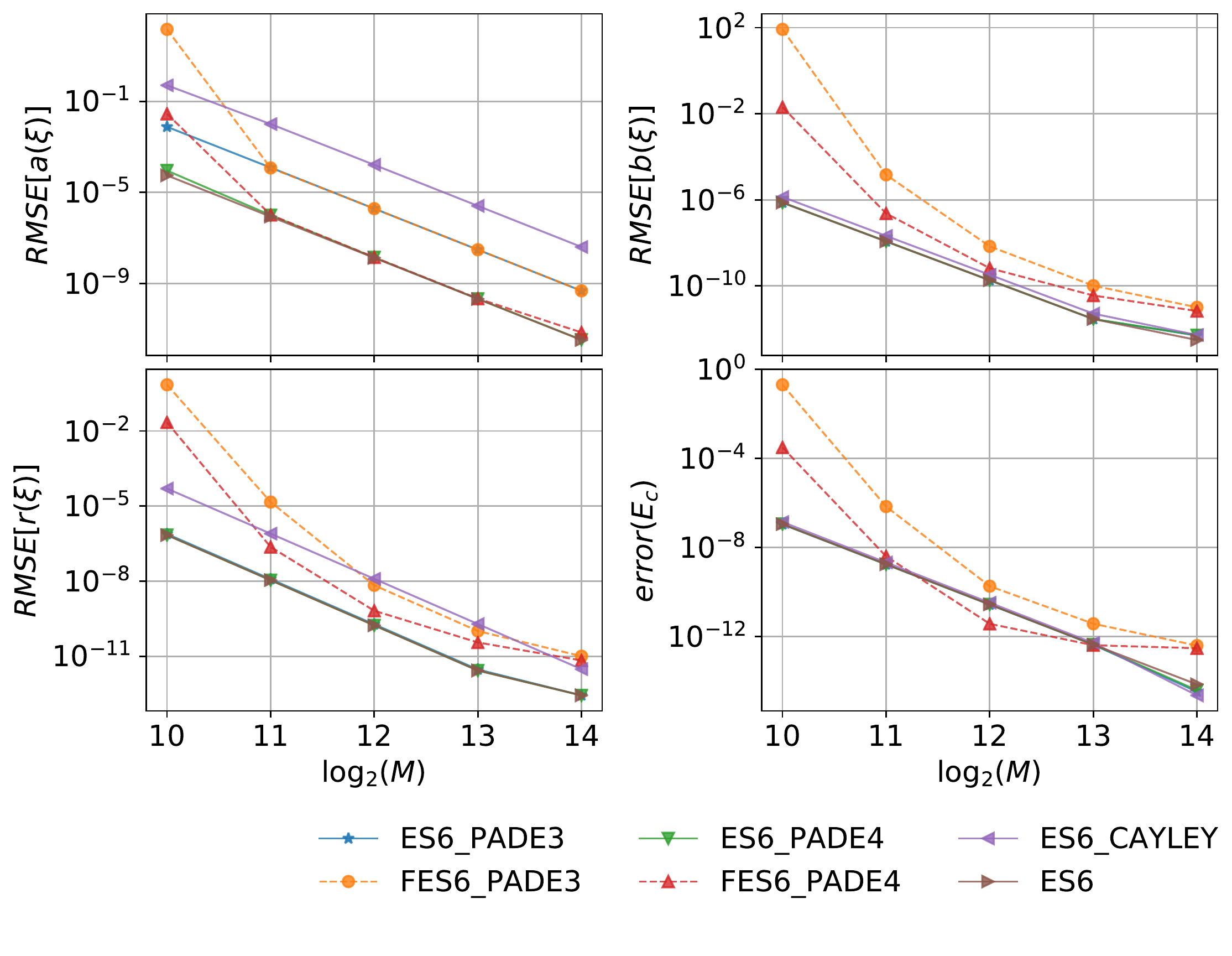}
	\caption{
		Continuous spectrum errors in the case of normal dispersion $\sigma=-1$. 
	}
	\label{fig:error2}
\end{figure}

Figures \ref{fig:error1} and \ref{fig:error2} present the continuous spectrum errors calculated using the schemes under consideration for the anomalous and normal dispersion, respectively. The best accuracy is shown by the schemes ES6 and ES6\_Pade4. ES6\_Pade3 is less accurate in calculating the coefficient $a(\xi)$. The worst result is obtained by the ES6\_Cayley. The fast schemes demonstrate the accuracy that is close to one of the initial schemes. 
Figures \ref{fig:error1} and \ref{fig:error2} also show the numerical errors for the continuous spectrum energy $E_c$ (\ref{Econt}). The accuracy of the fast schemes in calculating $E_c$ is worse, then the one of the conventional schemes. For normal dispersion the fast schemes are more accurate in computing $E_c$.

\begin{figure}[htpb]
	\centering
	\includegraphics[width=13cm]{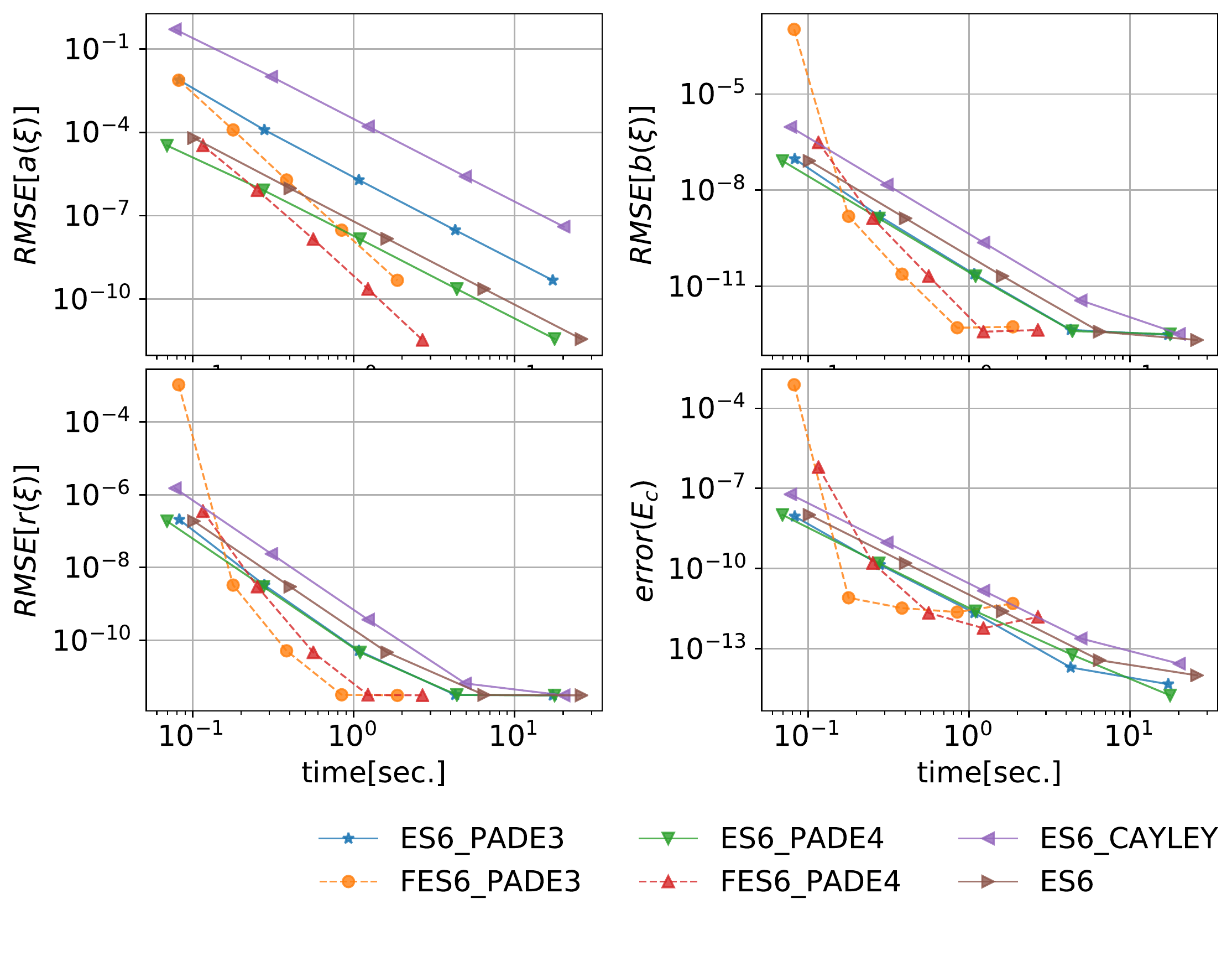}
	\caption{
		Continuous spectrum errors depending on the execution time trade-off in the case of anomalous dispersion $\sigma=1$.
	}
	\label{fig:errorVsTime1}
\end{figure}

\begin{figure}[htpb]
	\centering
	\includegraphics[width=13cm]{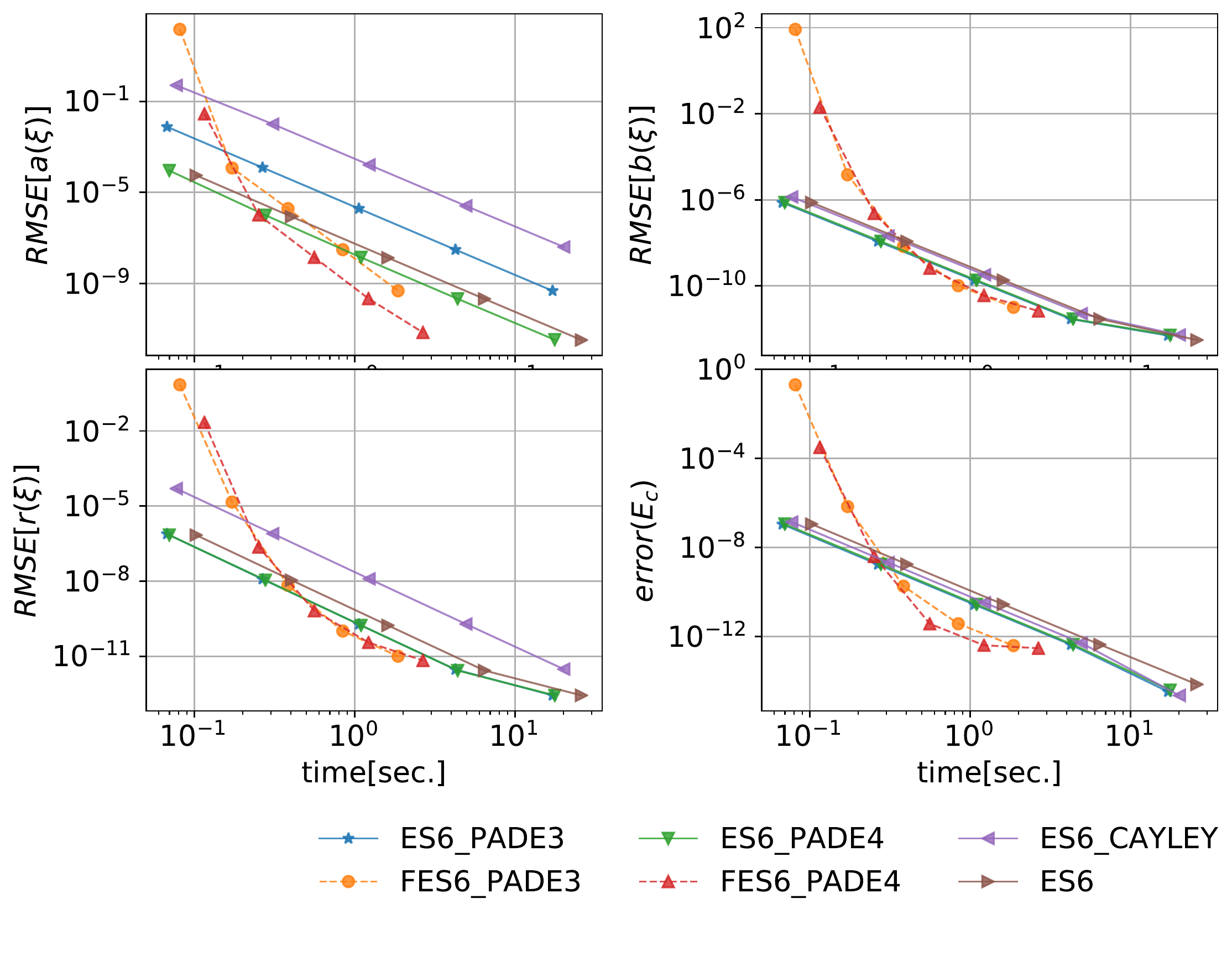}
	\caption{
			Continuous spectrum errors depending on the execution time trade-off in the case of normal dispersion $\sigma=-1$.
	}
	\label{fig:errorVsTime2}
\end{figure}

The efficiency of the schemes is compared in Figures \ref{fig:errorVsTime1} and \ref{fig:errorVsTime2}, where the continuous spectrum errors with respect to the running time are presented for the anomalous and normal dispersion, respectively. Among the conventional schemes the best result was obtained for the scheme with the fourth-order Pad\'e approximation ES6\_Pade4. The least efficient is the ES6\_Cayley scheme. 
The fast schemes outperform the conventional ones for a large number of nodes.
It is explained by the asymptotic complexity of the fast algorithms.
The FES6\_Pade3 has a smaller degree of the polynomial used for the transition matrix representation, so it works faster than the FES6\_Pade4. 
But the FES6\_Pade3 is less efficient in calculating $a(\xi)$ due to lack of accuracy.

\begin{figure}[htpb]
	\centering
	\includegraphics[width=13cm]{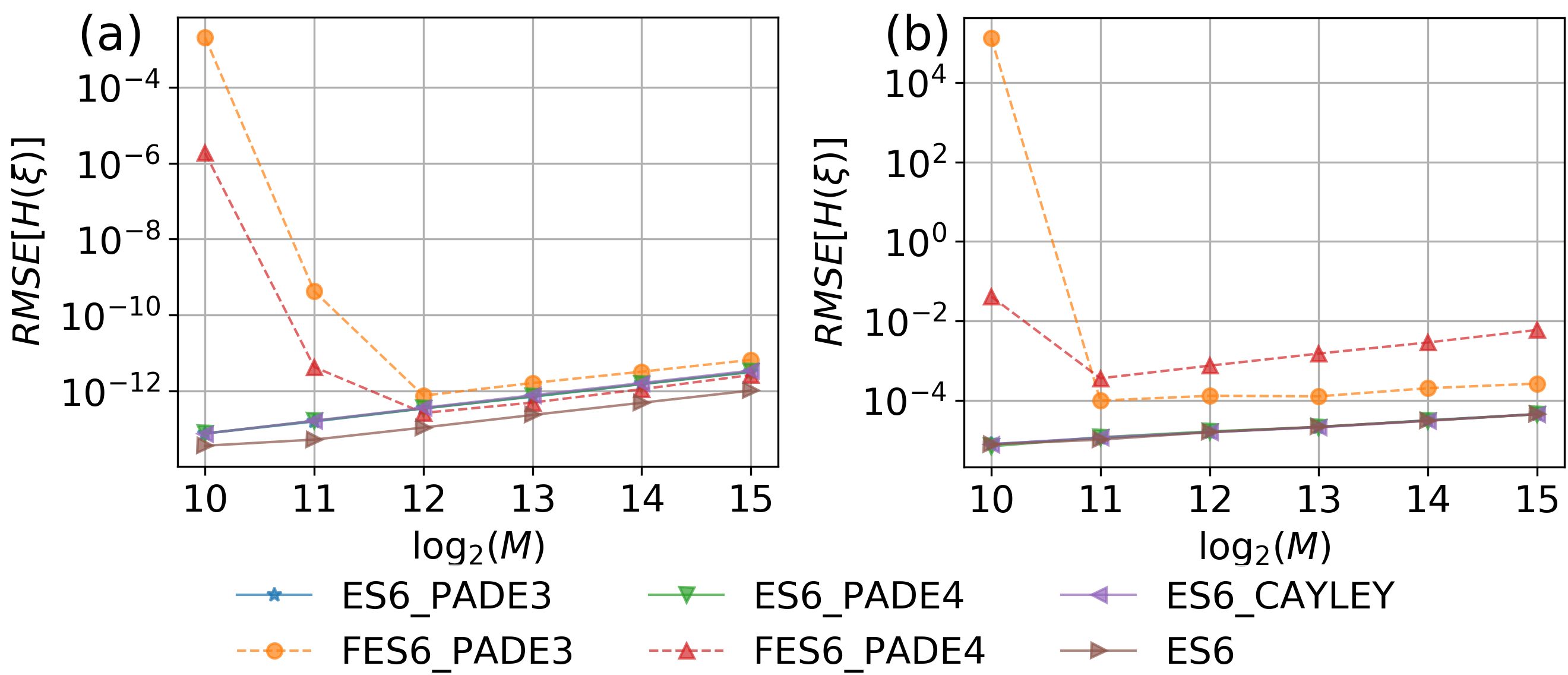}
	\caption{
		Invariant conservation error for anomalous dispersion $\sigma=1$ (a) and normal dispersion $\sigma=-1$ (b). 
	}
	\label{fig:rmseInvariant}
\end{figure}

\begin{figure}[htpb]
	\centering
	\includegraphics[width=0.8\textwidth]{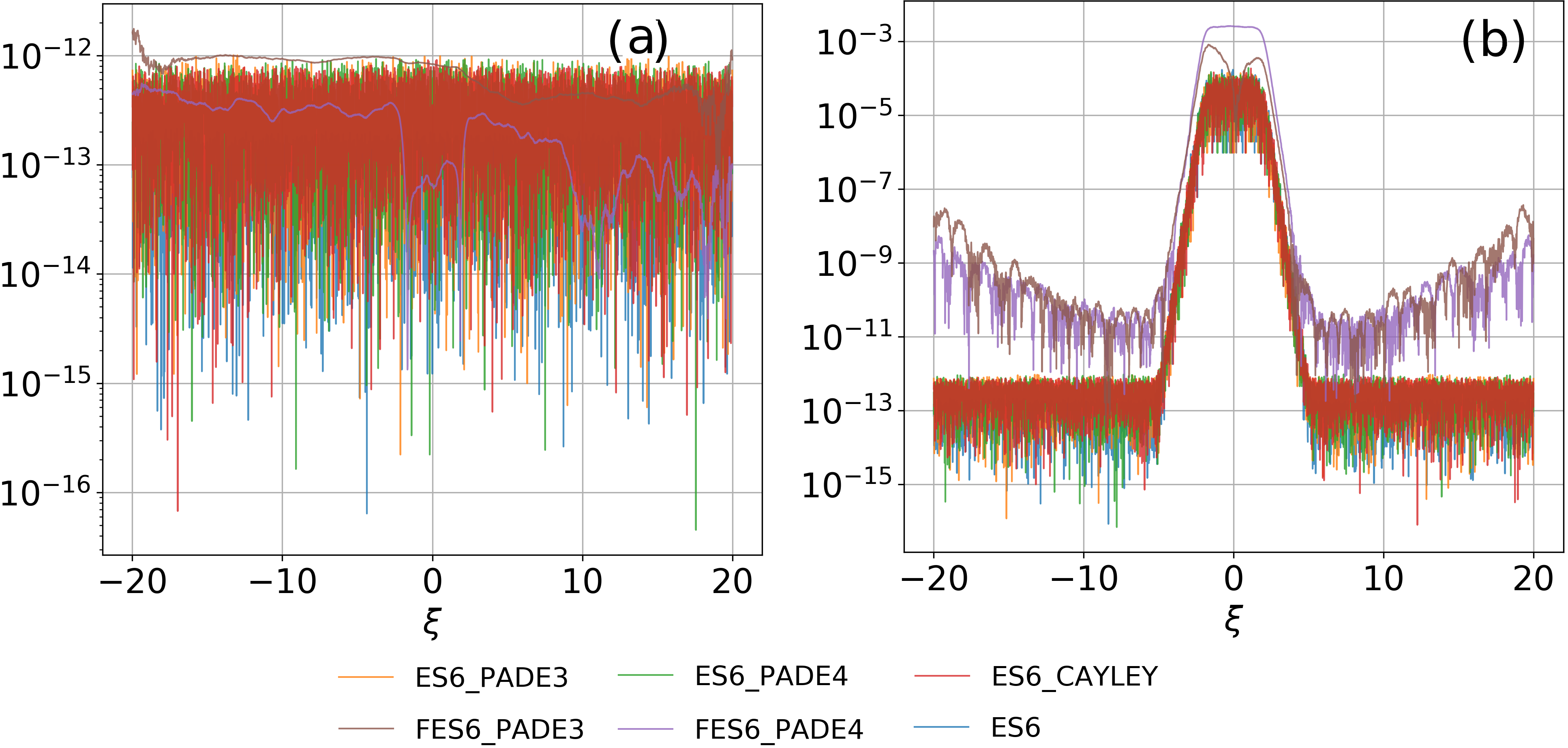}
	\caption{
		Invariant conservation error for anomalous dispersion $\sigma=1$ (a) and normal dispersion $\sigma=-1$ (b).
	}
	\label{fig:invariant}
\end{figure}

The conservation properties of the schemes are considered in Figures \ref{fig:rmseInvariant} and \ref{fig:invariant}. 
The quadratic invariant $H(\xi)$ is defined by (\ref{QuadInvariant}). All the conventional schemes demonstrate good conservation of the quadratic invariant.

Figure \ref{fig:rmseInvariant} presents the root mean squared error (\ref{MSE}) of $H(\xi)$ with respect to the number of points $M$ of the signal discretization. For anomalous dispersion, the fast schemes are close to the conventional ones starting from $M=2^{12}$. The FES6\_Pade4 is slightly better than the FES6\_Pade3. For normal dispersion, the fast schemes show worse results, and in this case, the FES6\_Pade3 works better than the FES6\_Pade4.

Figure \ref{fig:invariant} shows the error (\ref{err}) of calculating $H(\xi)$ with respect to the spectral parameter $\xi$ for $M=2^{12}$. In the case of anomalous dispersion the quadratic invariant $H(\xi)$ equally conserves for all  schemes considered here. In the case of normal dispersion the fast algorithms increase the error about one order of magnitude in the middle of the spectral interval and up to fourth order in the edges.
For normal dispersion, an error of all the schemes increases sufﬁciently in the middle of the spectral interval due to the subtraction of large modulo quantities.

\begin{figure}[htpb]
	\centering
	\includegraphics[width=13cm]{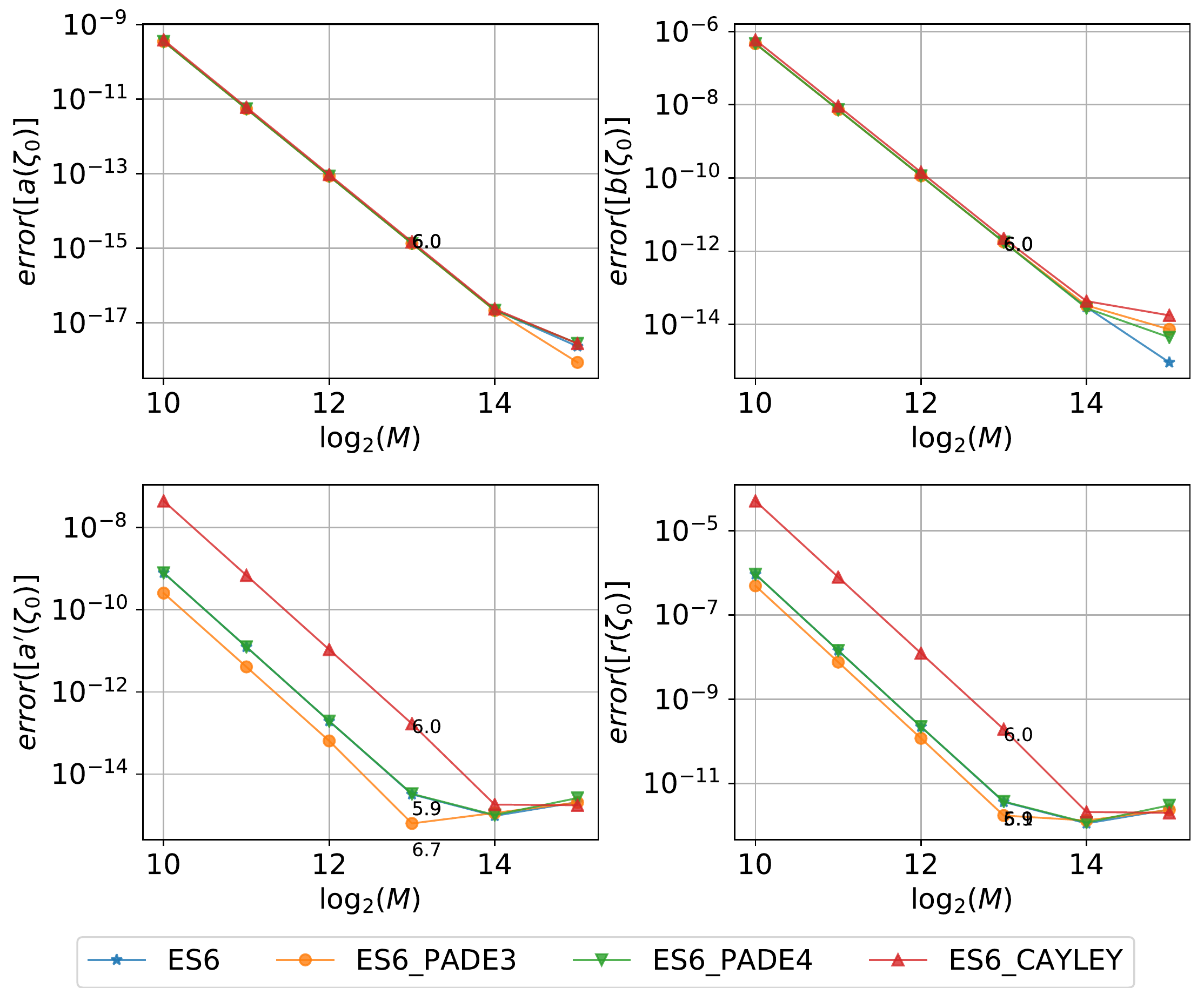}
	\caption{
		Discrete spectrum errors for the maximum eigenvalue $\zeta_0$.
	}
	\label{fig:discrete}
\end{figure}

The discrete spectrum errors are presented in Figure~\ref{fig:discrete}. Here we did not use any numerical algorithm for finding eigenvalues $\zeta_k$. The coefficients $a(\zeta_k)$, $b(\zeta_k)$, and $r(\zeta_k)$ were computed for the analytically known eigenvalues~\cite{medvedev2019exponential,medvedev2020exponential}. The review of the approaches for finding the eigenvalues can be found in recent papers~\cite{vasylchenkova2019direct, chimmalgi2019fast, chekhovskoy2020introducing-arxiv}.
Figure~\ref{fig:discrete} demonstrates the results calculated for the maximum eigenvalue $\zeta_0$. The coefficients $a(\zeta_0)$ and $b(\zeta_0)$ of the discrete spectrum are computed with almost the same accuracy for all the schemes. But for the derivative $a'(\zeta_0)$ and the phase coefficient $r_0$, the best result is obtained by the ES6\_Pade3 and the worst one by the ES6\_Cayley.

There are well-known problems with the computation of the coefficient $b(\zeta_k)$. We used the bi-directional algorithm~\cite{7486016} to find it. The algorithm is based on using
both boundary conditions~(\ref{psi0}) and~(\ref{psi0right}) to calculate the coefficient $b(\zeta_{k} )$ of the discrete spectrum:
\begin{equation}\label{bidir}
 \Psi(t,\zeta_k)=\Phi(t,\zeta_k)b(\zeta_k).
\end{equation}

To find the phase coefficients $r(\zeta_k)$ we need to know the derivative of the transition matrix $T'$
with respect to the spectral parameter $\zeta$. The formulas can be found in Appendix \ref{appendix_derivative}.

\section{Conclusion}

Families of schemes of the sixth order are constructed for a system of linear differential equations of the first order with a matrix depending on time and spectral parameter. Such schemes are supposed to be used in the numerical solution of the direct spectral problem for integrable vector nonlinear Schr\"odinger equations; therefore, the main attention was paid to schemes that allow the use of fast algorithms when solving the system for a large number of spectral parameter values. In particular, the proposed schemes are applied to solve the direct spectral problem for the ZS system. In our opinion, the constructed schemes will be useful for more accurate realistic calculations in the construction of telecommunication data transmission systems based on NLSE soliton solutions. On the other hand, the proposed schemes of the sixth order of accuracy are on the verge of computational consistency, because the schemes of the next order of accuracy require a large number of points for approximating the derivatives and contain a significantly larger number of terms.

\section*{Funding}
Russian Science Foundation (RSF) (20-11-20040). 


\appendix


\section{Derivative of the transition matrix}\label{appendix_derivative}
To find the phase coefficients $r(\zeta_k)$ we need to know the derivative $a'(\zeta)$ at the point $\zeta=\zeta_k$
\begin{equation}
\frac{da}{d\zeta} = \frac{d\psi_1}{d\zeta}e^{i\zeta (L-\tau/2)} + i(L-\tau/2) a(\zeta).
\end{equation}
From (\ref{T0}) we get
\begin{equation}\label{dpsi}
\frac{d}{d\zeta}\Psi_{n+\frac{1}{2}} = T' \Psi_{n-\frac{1}{2}} + T\frac{d}{d\zeta}\Psi_{n-\frac{1}{2}},
\end{equation}
where the initial value is defined from (\ref{psi0})
\begin{equation}\label{dpsi0}
\frac{d}{d\zeta}\Psi(-L-\tau/2,\zeta)=\left(\begin{array}{c}-i (-L-\tau/2)\psi_1(-L-\tau/2,\zeta)\\0\end{array}\right)
\end{equation}
and the derivative $T'$ with respect to the spectral parameter $\zeta$ has the form
\begin{equation}
T'=c'\sigma_0+\left(\frac{s(\lambda)}{\lambda}\right)'Z+\frac{s(\lambda)}{\lambda}Z',
\end{equation}
where
\begin{equation}
c'=\frac{dc(\lambda)}{d\lambda}\,\frac{z_1z_1'+z_2z_2'+z_3z_3'}{\lambda},
\end{equation}
\begin{equation}
\left(\frac{s(\lambda)}{\lambda}\right)'=\left(\frac{ds(\lambda)}{d\lambda}-\frac{s(\lambda)}{\lambda}\right)\frac{z_1z_1'+z_2z_2'+z_3z_3'}{\lambda^2}.
\end{equation}

The derivatives $dc_k/d\lambda$ and $ds_k/d\lambda$ for Pad\'e approximations $E_k(Z)$ are easy to find.

\bibliographystyle{unsrt} 
\bibliography{magnus,sample,references} 

\end{document}